\documentclass[twoside,reqno,A4]{amsart}

\usepackage{latexsym,rotate,eucal,cite,url}
\usepackage{amsmath,amsthm,amssymb,amsxtra} 

\newcommand{\quab}{\hspace*{2.0mm}}
\newcommand{\puab}{\hspace*{-1.5mm}}

\newcommand{\sst}{\scriptstyle}

\newcommand{\qqed}{\tag*{\qed}} 

\newcommand{\pav}[1]{\lfloor{#1}\rfloor}

\newcommand{\pavv}[1]{\left\lfloor{#1}\right\rfloor}

\newcommand{\ang}[1]{\langle{#1}\rangle}

\newcommand{\mb}[1]{\mathbb{#1}}
\newcommand{\mc}[1]{\mathcal{#1}}

\newcommand{\mv}{\mathversion}

\newcommand{\hh}[2]{\left[#1\right]_{#2}}
\newcommand{\hyp}[3]{\left[\puab\ba{#1}#2\\#3\ea\puab\right]}
\newcommand{\hyq}[4]{\left[\puab\ba{#1}#3\\#4\ea{\puab\Big|#2}\right]}

\newcommand{\binm}{\binom}

\newcommand{\nnm}{\nonumber}
\newcommand{\be}{\begin{equation}}
\newcommand{\ee}{\end{equation}}
\newcommand{\ba}{\begin{array}}
\newcommand{\ea}{\end{array}}
\newcommand{\bmn}{\begin{eqnarray}}
\newcommand{\emn}{\end{eqnarray}}
\newcommand{\bnm}{\begin{eqnarray*}}
\newcommand{\enm}{\end{eqnarray*}}
\newcommand{\bln}{\begin{subequations}}
\newcommand{\eln}{\end{subequations}}

\newcommand{\pq}[1]{\begin{equation}#1\end{equation}}
\newcommand{\pp}[2]{\begin{aligned}[#1]#2 
            \end{aligned}}  

\newcommand{\pmq}[1]{\begin{align}#1
            \end{align}}     
\newcommand{\pnq}[1]{\begin{align*}#1
            \end{align*}}    
\newcommand{\pnp}[2]{\begin{alignat*}{#1}#2
            \end{alignat*}}  
\newcommand{\centro}[1]
           {\begin{center}#1\end{center}}

\newcommand{\alp}{\alpha}
\newcommand{\bet}{\beta}
\newcommand{\gam}{\gamma}

\newcommand{\lam}{\lambda}

\newcommand{\veps}{\varepsilon}

\newcommand{\sig}{\sigma}

\newcommand{\Del}{\Delta}

\newcommand{\Lam}{\Lambda}
\newcommand{\Ome}{\Omega}

\newcommand{\Gam}{\Gamma}

\newtheorem{thm}{Theorem}
\newtheorem{lemm}[thm]{Lemma}

\newtheorem{prop}[thm]{Proposition}
\newtheorem{exam}{Example}

\newtheorem{entry}{Entry}

\newcommand{\referxy}[4]{\bibitem{kn:#1}{#2,}~\emph{#3,}~{#4.}}	
\newcommand{\cito}[1]{\cite{kn:#1}}	
\newcommand{\citu}[2]{\cite[#2]{kn:#1}}

\begin{document}{\hfill\fbox{\textbf{First Version}}} 
\title{Triplicate Dual Series of Dougall--Dixon Theorem}
\author{Xiaojing Chen and Wenchang Chu}
\address{School of Statistics\newline
	 Qufu Normal University\newline
	 Qufu (Shandong), P.~R.~China}
\email{upcxjchen@outlook.com}
\address{Department of Mathematics and Physics\newline
     	 University of Salento (P.~O.~Box~193) \newline
      	 73100 Lecce, ~ Italy}
\email{chu.wenchang@unisalento.it}
\thanks{Email addresses: upcxjchen@outlook.com
	and chu.wenchang@unisalento.it}


\begin{abstract}
Applying the triplicate form of the extended Gould--Hsu inverse
series relations to Dougall's summation theorem for the well--poised
$_7F_6$-series, we establish, from the dual series, several interesting
Ramanujan--like infinite series expressions for $\pi^2$ and $\pi^{\pm1}$
with convergence rate ``$-\frac{1}{27}$".
\end{abstract}

\maketitle\thispagestyle{empty}
\markboth{Xiaojing Chen and Wenchang Chu}
{Triplicate Dual Series of Dougall--Dixon Theorem} 

\vspace*{-5mm}\section{Introduction and Motivation}
For an indeterminate $x$, the shifted factorial is defined by
$(x)_0\equiv1$ and
\[(x)_n=x(x+1)\cdots(x+n-1)
\quad\text{for}\quad
n\in\mb{N}.\]
This can also be expressed as a quotient
$(x)_n=\Gam(x+n)/\Gam(x)$, where the $\Gam$-function
(see~\citu{rain}{\S8} for example) is given by the
beta integral
\[\Gam(x)=\int_{0}^{\infty}u^{x-1}
e^{-u}\mathrm{d}u
\quad\text{for}\quad
\Re(x)>0,\]
which admits Euler's reflection property
\pq{\Gam(x)\Gam(1-x)=\label{reflex}
\frac{\pi}{\sin\pi x}
\quad\text{with}\quad \Gam(\tfrac12)=\sqrt{\pi}}
and the following asymptotic formula
\pq{\Gam(x+n)\approx n^{x}(n-1)!
\quad\text{as}\quad\label{limit}
n\to\infty.}
This last formula is simpler than Stirling formula
and utilized frequently to evaluate limits of
$\Gam$-function quotients.

About one century ago, Dougall~\cito{dougall} discovered a very general
summation theorem for the terminating well--poised $_7F_6$-series.
By making use of the triplicate form of the extended Gould--Hsu inverse
series relations, we shall investigate the dual series of Dougall's
well--poised sum, that will lead to a large class of summation formulae
for $\pi$-related infinite series of convergence rate ``$\frac1{729}$".
According to the bisection series method, a number of these series
can be reduced to simpler ones with convergence rate ``$-\frac1{27}$".
Five elegant formulae are anticipated as follows:
\pnq{\tag{Example~\ref{pi-1A}}
\frac{9\sqrt3}{\pi\sqrt[3]{16}}
&=\sum_{k=0}^{\infty}
\Big(\frac{-1}{27}\Big)^k
\hyp{cccccc}{\frac13,\:\frac23,\:\frac16}
{1,\:1,\:1\rule[2mm]{0mm}{2mm}}_k
\big\{2 + 21k\big\},\\
\tag{Example~\ref{pi-1B}}
\frac{27\sqrt3}{\pi\sqrt[3]{32}}
&=\sum_{k=0}^{\infty}
\Big(\frac{-1}{27}\Big)^k
\hyp{cccccc}{\frac13,\:\frac23,\:\frac56}
{1,\:1,\:1\rule[2mm]{0mm}{2mm}}_k
\big\{5 + 42k\big\},\\
\tag{Example~\ref{pi-1C}}
\frac{3\sqrt3}{\pi\sqrt[3]4}
&=\sum_{k=0}^{\infty}
\Big(\frac{-1}{27}\Big)^k
\hyp{cccccc}{\frac13,-\frac13,-\frac16}
{1,~1,~1\rule[2mm]{0mm}{2mm}}_k
\big\{1-63 k^2\big\},\\
\tag{Example~\ref{pi=4/3}}
\frac{\pi\Gamma^2 (\frac{1}{3})}{6\Gamma^2(\frac{5}{6})}
&=\sum_{k=0}^{\infty}
\Big(\frac{-1}{27}\Big)^k
\hyp{cccccc}{1,\,\frac{1}2,\frac{2}3}
{\frac43,\frac43,\frac43\rule[2mm]{0mm}{2mm}}_k
\big\{3+7 k\big\},\\
\tag{Example~\ref{pi=5/3}}
\frac{48\pi\Gamma^2 (\frac{2}{3}) }{\Gamma^2(\frac{1}{6})}
&=\sum_{k=0}^{\infty}
\Big(\frac{-1}{27}\Big)^k
\hyp{cccccc}{1,\,\frac{1}2,\frac{1}3}
{\frac53,\frac{5}3,\frac{5}3\rule[2mm]{0mm}{2mm}}_k
\big\{9+28 k+21 k^2\big\}.}
They resemble the so--called Ramanujan--like series, mainly
discovered one century ago by Ramanujan~\cito{ramanujan} and
recently by Guillera~\cite{kn:gj03em,kn:gj06rj,kn:gj08rj}.

The rest of the paper will be organized as follows.
The next section will serve as the theoretical part, where
the main theorems and proofs will be included. Then in
Section~3, we shall present 35 infinite series expressions
for $\pi^2$ and $\pi^{\pm1}$ as applications.

Throughout the paper, the following abbreviated notations
will be adopted for product and quotient of shifted factorials:
\pnq{[\alp,\bet,\cdots,\gam]_n\quab
=&(\alp)_n(\bet)_n\cdots(\gam)_n,\\
\hyp{cccc}{\alp,\bet,\cdots,\gam}
          {A,B,\cdots,C}_n
=&\frac{(\alp)_n(\bet)_n\cdots(\gam)_n}
        {(A)_n(B)_n\cdots(C)_n}.}

\section{Triplicate Inversion of Dougall's \mv{bold}$_7F_6$-Series}
In 1973, Gould and Hsu~\cito{hsu} discovered a useful pair
of inverse series relations, which can equivalently be
reproduced below. Let $\{a_i,b_i\}$ be any two complex
sequences such that the $\phi$-polynomials defined by
\pq{\phi(x;0)\equiv1\label{phi-pp}
\quad\text{and}\quad
\phi(x;n)=\prod_{k=0}^{n-1}(a_k+xb_k)
\quad\text{for}\quad n\in\mb{N}}
differ from zero for $x,\:n\in\mb{N}_0$.
Then there hold the inverse series relations
\pmq{
f(n)&=\label{inv+gh}
\sum_{k=0}^{n}
(-1)^k\binm{n}{k}
\phi(k;n)~g(k),\\
g(n)&=\label{inv-gh}
\sum_{k=0}^{n}
(-1)^k\binm{n}{k}
\frac{a_k+kb_k}{\phi(n;k+1)}
f(k).}

The Gould--Hsu inversions have the following extended form
(cf.~\cite{kn:bress,kn:chu93b,kn:chu94e}):
\pmq{
f(n)&=\label{inv+w}
\sum_{k=0}^{n}(-1)^k\binm{n}{k}
\phi(\lam+k;n)\phi(-k;n)
\frac{\lam+2k}{(\lam+n)_{k+1}}g(k),\\
g(n)&=\label{inv-w}
\sum_{k=0}^{n}
(-1)^k\binm{n}{k}
\frac{(a_k+\lam b_k+kb_k)(a_k-kb_k)}{\phi(\lam+n;k+1)\phi(-n;k+1)}(\lam+k)_n
f(k);}
where the `adjunct' factor can be expressed transparently as
\[(a_k+\lam b_k+kb_k)(a_k-kb_k)
=\frac{\phi(\lam+k;k+1)\phi(-k;k+1)}{\phi(\lam+k;k)\phi(-k;k)}.\]

These inverse series relations have been shown powerful 
in dealing with terminating hypergeometric series identities
\cite{kn:chu93b,kn:chu93c,kn:chu94a,kn:chu94b,kn:chu94e,kn:gess82s}.
Their duplicate forms and the corresponding $q$-analogues due to 
Carlitz~\cito{carlitz} with respective applications were extensively 
explored in \cite{kn:chu13e,kn:chu00a,kn:chu02b,kn:chu08j} and 
\cite{kn:chu94c,kn:chu94d,kn:chu95b,kn:chu06d,kn:gess83s,kn:gess86s}.

By employing the above inverse pair, we shall work out several new 
$\pi$-related infinite series expressions. Recall the fundamental 
identity discovered by Dougall~\cito{dougall}
(see also Bailey~\citu{bailey}{\S4.3}) for very--well--poised
terminating $_7F_6$-series
\pq{\label{doug}
\pp{c}{\Ome_n(a;b,c,d)
~:=~\hyp{c}{1+a,1+a-b-c,1+a-b-d,1+a-c-d}{1+a-b,1+a-c,1+a-d,1+a-b-c-d}_n&\\
=\sum_{k=0}^n\frac{a+2k}{a}\hyp{ccccc}
{a,~b,~c,~d,~e,\:-n}
{1,1+a-b,1+a-c,1+a-d,1+a-e,1+a+n}_k&,}}
where the series is 2-balanced because $1+2a+n=b+c+d+e$.

For all $n\in\mb{N}_0$, it is well known that
$\boxed{n=\pavv{\tfrac{n}3}+\pavv{\tfrac{1+n}3}+\pavv{\tfrac{2+n}3}}$, where $\pav{x}$
denotes the greatest integer not exceeding $x$.
Then it is not difficult to check that Dougall's formula \eqref{doug}
is equivalent to the following one

\pnq{&\Ome_n\big(a;b+\pavv{\tfrac{n}3},c+\pavv{\tfrac{1+n}3},d+\pavv{\tfrac{2+n}3}\big)
=\frac{(1+a)_n}{(b+c+d-a)_n}\\
&\times\frac{(1+a-c-d)_{\pavv{\frac{n}3}}}{(b-a)_{\pavv{\frac{n}3}}}
\frac{(1+a-b-d)_{\pavv{\frac{1+n}3}}}{(c-a)_{\pavv{\frac{1+n}3}}}
\frac{(1+a-b-c)_{\pavv{\frac{2+n}3}}}{(d-a)_{\pavv{\frac{2+n}3}}}\\
&\times\frac{(b+c-a)_{\pavv{\frac{2n}3}}}{(1+a-d)_{\pavv{\frac{2n}3}}}
\frac{(b+d-a)_{\pavv{\frac{1+2n}3}}}{(1+a-c)_{\pavv{\frac{1+2n}3}}}
\frac{(c+d-a)_{\pavv{\frac{2+2n}3}}}{(1+a-b)_{\pavv{\frac{2+2n}3}}},}
with its parameters subject to $\boxed{1+2a=b+c+d+e}$.

Reformulate the above equality as a binomial sum
\pnq{\sum_{k=0}^n&(-1)^k\binm{n}{k}
\frac{a+2k}{(a+n)_{k+1}}
\hyp{c}{a,~b,~c,~d,~1+2a-b-c-d}{1+a-b,1+a-c,1+a-d,b+c+d-a}_k\\
&\times\hh{b+k,b-a-k}{\pavv{\frac{n}3}}
\hh{c+k,c-a-k}{\pavv{\frac{1+n}3}}
\hh{d+k,d-a-k}{\pavv{\frac{2+n}3}}\\
=~&\hh{b,1+a-c-d}{\pavv{\frac{n}3}}
\hh{c,1+a-b-d}{\pavv{\frac{1+n}3}}
\hh{d,1+a-b-c}{\pavv{\frac{2+n}3}}\\
&\times\frac{(a)_n}{(b+c+d-a)_n}
\hyp{c}{b+c-a}{1+a-d}_{\pavv{\frac{2n}3}}
\hyp{c}{b+d-a}{1+a-c}_{\pavv{\frac{1+2n}3}}
\hyp{c}{c+d-a}{1+a-b}_{\pavv{\frac{2+2n}3}}.}

This equality matches exactly to \eqref{inv+w} under the assignments
$\lam\to a$ and
\[\phi(x;n)\to(b-a+x)_{\pavv{\frac{n}3}}(c-a+x)_{\pavv{\frac{1+n}3}}
            (d-a+x)_{\pavv{\frac{2+n}3}}\]
as well as
\pnq{f(n)&\to n!~(a)_n\times\mc{F}(n),\\
g(k)&\to\hyp{c}{a,~b,~c,~d,~1+2a-b-c-d}{1+a-b,1+a-c,1+a-d,b+c+d-a}_k;}
where
\pmq{\label{mc-ff}\mc{F}(n)&=\hyp{c}
{b+c-a}{1+a-d}_{\pav{\frac{2n}3}}
\hyp{c}{b+d-a}{1+a-c}_{\pav{\frac{1+2n}3}}
\hyp{c}{c+d-a}{1+a-b}_{\pav{\frac{2+2n}3}}\\
&\nnm\times\frac{\hh{b,1+a-c-d}{\pav{\frac{n}3}}
\hh{c,1+a-b-d}{\pav{\frac{1+n}3}}
\hh{d,1+a-b-c}{\pav{\frac{2+n}3}}}{n!~(b+c+d-a)_n}.}

For the sake of brevity, we introduce the $\psi$-polynomials by
\pq{\label{psi-pp}\pp{c}{\psi(x;n)=\phi(a+x;n)\phi(-x;n)
=(b+x)_{\pav{\frac{n}3}}(c+x)_{\pav{\frac{1+n}3}}(d+x)_{\pav{\frac{2+n}3}}&\\
\times(b-a-x)_{\pav{\frac{n}3}}(c-a-x)_{\pav{\frac{1+n}3}}(d-a-x)_{\pav{\frac{2+n}3}}&.}}

Then the dual relation corresponding to \eqref{inv-w}
can explicitly be stated, after some simplifications,
in the following lemma.
\begin{lemm}\label{tre-A}
For the $\mc{F}$-quotient of shifted factorials and the $\psi$-polynomials
defined respectively in \eqref{mc-ff} and \eqref{psi-pp}, we have the
summation formula
\pnq{&\hyp{c}{b,~c,~d,~1+2a-b-c-d}
{1+a-b,1+a-c,1+a-d,b+c+d-a}_n\\
&=\sum_{k=0}^n\mc{F}(k)
\frac{\psi(k;k+1)}{\psi(k;k)}
\frac{(-n)_k(a+n)_k}{\psi(n;k+1)}.}
\end{lemm}

Observe that $\psi(n;k+1)$ is a polynomial of degree $2k+2$ in $n$ with
the leading coefficient equal to $(-1)^{k+1}$.
Now multiply by $``n^2"$ across the binomial relation in Lemma~\ref{tre-A}
and then let $n\to\infty$. We may evaluate the limits of the left member
by \eqref{limit} and of the corresponding right member through Weierstrass's
$M$-test on uniformly convergent series (cf. Stromberg~\citu{karl}{\S3.106}).
After some routine simplifications, the resultant limiting relation can be
expressed explicitly as follows.
\begin{prop}\label{tre-B}
Let $\Gam(a,b,c,d)$ stand for the quotient of the $\Gam$-function given by
\pq{\label{gamma}\Gam(a,b,c,d)
=\frac{\Gam(1+a-b)\Gam(1+a-c)\Gam(1+a-d)\Gam(b+c+d-a)}
      {\Gam(b)~\Gam(c)~\Gam(d)~\Gam(1+2a-b-c-d)}.}
Then for the $\mc{F}$-quotient of shifted factorials and the $\psi$-polynomials
defined respectively in \eqref{mc-ff} and \eqref{psi-pp}, the following
infinite series identity holds:
\[\Gam(a,b,c,d)
=-\sum_{k=0}^{\infty}
\frac{\psi(k;k+1)}{\psi(k;k)}
\mc{F}(k).\]
\end{prop}

Let $\veps=0,1,2$ and $T(k)$ be the summand with index $k$ in the last series.
Putting the initial $\veps$ terms aside and then classifying the remaining
terms with respect to their indices modulo $3$, we get the expressions
\pnq{\sum_{k=0}^{\infty}T(k)
&=\sum_{k=0}^{\veps-1}T(k)
+\sum_{k=0}^{\infty}\sum_{i=0}^2T(\veps+i+3k)\\
&=\sum_{k=0}^{\veps-1}T(k)
+\sum_{k=1}^{\infty}\sum_{j=1}^3T(\veps-j+3k).}

Denote further by $\sig(\veps)$, $\Del_k(\veps)$ and $\nabla_k(\veps)$ the sum of initial
$\veps$-terms and the ``weight functions" (where the latter are clearly a rational
function of $k$):
\pmq{\sig(\veps)~&=\label{sigma}\sum_{k=0}^{\veps-1}
\bigg\{\frac{\psi(k;k+1)}{-\psi(k;k)}\bigg\}\mc{F}(k),\\
\Del_k(\veps)&=\label{delta}\sum_{i=0}^2
\bigg\{\frac{\psi(\veps+i+3k;1+\veps+i+3k)}
{-\psi(\veps+i+3k;\veps+i+3k)}\bigg\}
\frac{\mc{F}(\veps+i+3k)}{\mc{F}(3k)},\\
\nabla_k(\veps)&=\label{nabla}\sum_{j=1}^3
\bigg\{\frac{\psi(\veps-j+3k;1+\veps-j+3k)}
{-\psi(\veps-j+3k;\veps-j+3k)}\bigg\}
\frac{\mc{F}(\veps-j+3k)}{\mc{F}(3k)}.}
Then the identity in Proposition~\ref{tre-B} can be restated
in the theorem below.
\begin{thm}\label{tre-C}
Assume that $\Gam(a,b,c,d)$, $\sig(\veps)$, $\Del_k(\veps)$ and $\nabla_k(\veps)$
are as in \eqref{gamma}, \eqref{sigma}, \eqref{delta} and \eqref{nabla} respectively.
Then for $\veps=0,1,2$ and the $\mc{F}$-quotient of shifted factorials defined
in \eqref{mc-ff}, the following infinite series identities hold:
\pnq{\Gam(a,b,c,d)
&=\sig(\veps)+\sum_{k=0}^{\infty}
\Del_k(\veps)\mc{F}(3k)\\
&=\sig(\veps)+\sum_{k=1}^{\infty}
\nabla_k(\veps)\mc{F}(3k).}
\end{thm}

In the above theorem, the series is expressed in two different manners
because it happens frequently that a series with its summation index
initiating at $k=0$ has better looking than that at $k=1$, or vice
versa. This will be seen from our examples in the next section.
In the above series, $\mc{F}(3k)$ results in the dominant part
\pnq{\mc{F}(3k)&=\hyp{c}
{b+c-a,b+d-a,c+d-a}
{1+a-b,1+a-c,1+a-d}_{2k}\\
&\times\frac{\hh{b,c,d,1+a-b-c,1+a-b-d,1+a-c-d}{k}}
{(3k)!~(b+c+d-a)_{3k}},}
which determines the convergence rate of the series
to be ``$\frac1{729}$". Instead, both $\Del_k(\veps)$ and $\nabla_k(\veps)$
are perturbing parts consisting of
only a few terms. Therefore for specific values of
$\veps$ and $\{a,b,c,d\}$, in order to find the
infinite series identity, it is enough to compute
the corresponding $\sig(\veps)$ and $\Del_k(\veps)$
(or $\nabla_k(\veps)$), and then to simplify
the resultant expression.

\section{Infinite Series of Ramanujan Type Involving \mv{bold}$\pi$}
By specifying the parameters $\{a,b,c,d\}$, we can derive numerous
infinite series identities of convergence rate ``$\tfrac1{729}$" from
Theorem~\ref{tre-C} with $\veps=0,1,2$. However in general, there
are complicated weight polynomials appearing in the summands
of these series. Two examples are illustrated as follows.
Letting
\[\veps=0\quad\text{and}\quad
\big\{a,b,c,d\big\}=\bigg\{\frac12,\frac12,\frac12,\frac14\bigg\},\]
we can compute with \emph{Mathematica} commands
\pnq{&\mc{F}(3k)=\hyp{c}
{\frac12,\frac14,\frac14}
{\\[-3.5mm]1,1,\frac54}_{2k}
\frac{\hh{\frac12,\frac12,\frac12,\frac14,\frac34,\frac34}{k}}
{(3k)!~(\frac34)_{3k}}
=\hyp{cccccc}{\frac{1}{2},\frac{1}{4},\frac{3}{4},\frac{3}{4},\frac{3}{4},~\frac{1}{8},\frac{1}{8},\frac{5}{8}}
{1,1,1,\frac{1}{3},\frac{2}{3},\frac98,\frac{7}{12},\frac{11}{12}\rule[2mm]{0mm}{2mm}}_k
\Big(\frac1{729}\Big)^k,\\
&\Del_k(0)=\frac{(1+8k)(93184
   k^4+154432 k^3+85840 k^2+17484 k+855)}{768(1+3k) (2+3k)(7+12k)},\\
&\sig(0)=0
\quad\text{and}\quad
\Gam(\tfrac12,\tfrac12,\tfrac12,\tfrac14)=\frac{1}{4\pi}.}
Substituting them into Theorem~\ref{tre-C} and then multiplying by ``$4\times2688$"
across the resultant equation, we derive the following infinite series identity.
\begin{exam}[\fbox{$\veps=0:\frac12,\frac12,\frac12,\frac14$}]
\[ \frac{2688}{\pi }=\sum_{k=0}^{\infty}
\hyp{cccccc}{\frac12,\frac14,\frac34,\frac34,\frac34,~\frac18,\frac58}
{1,1,1,\frac43,\frac53,\frac{11}{12},\frac{19}{12} \rule[2mm]{0mm}{2mm}}_k
\frac{\sst 93184 k^4+154432 k^3+85840 k^2+17484 k+855}{729^k}.\]
\end{exam}

Analogously, we have another infinite series identity of similar type.
\begin{exam}[\fbox{$\veps=1:\frac32,1,1,\frac56$}]
\[\frac{20\pi}{3\sqrt{3}}=10+\sum_{k=1}^{\infty}
\hyp{cccccc}{1,-\frac{1}{2},\frac{2}{3},-\frac{1}{3},-\frac{1}{3},\frac{1}{6},-\frac{5}{6}}
{\frac13,~\frac43,~\frac14,~\frac34,~\frac{4}{9},~\frac{7}{9},~\frac{10}{9} \rule[2mm]{0mm}{2mm}}_k
\frac{\sst 19656 k^4-7749 k^3-3150 k^2+613 k+118}{729^k}.\]
\end{exam}

When $\boxed{b+c+d-a}$ equals to a half integer, the corresponding series
in Theorem~\ref{tre-C} can be reformulated, by means of the bisection series
method, as a simpler series with convergence rate ``$\tfrac{-1}{27}$". In order
to show how this approach works, we present demonstrations in details
for two infinite series identities.

We start with the following strange evaluation of a hypergeometric $_3F_2$-series.
\begin{exam}[\fbox{$\veps=0:\frac53,\frac53,\frac23,\frac56$}]\label{3F2-S}
\[{_3F_2}\hyq{r}{\frac{-1}{27}}
{\frac23,\frac73,-\frac16}
{1,~2~\rule[2mm]{0mm}{-2mm}}
:=\sum_{k=0}^{\infty}
\hyp{ccc}{\frac23,\frac73,-\frac16}
{1,~1,~2\rule[-2mm]{0mm}{2mm}}_k
\Big(\frac{-1}{27}\Big)^k
=\frac{81\sqrt3}{28\cdot2^{2/3}\pi}.\]
\end{exam}

\emph{Proof.} \
By specifying the parameters in Theorem~\ref{tre-C}
\[\veps=0\quad\text{and}\quad
\big\{a,b,c,d\big\}=\bigg\{\frac53,\frac53,\frac23,\frac56\bigg\}\]
we have
\pnq{&\mc{F}(3k)=\hyp{c}
{\frac23,\frac56,-\frac16}
{\\[-3.5mm]1,2,~\frac{11}6}_{2k}
\frac{\hh{\frac13,\frac23,\frac53,\frac16,\frac56,\frac76}{k}}
{(3k)!~(\frac32)_{3k}}
=\hyp{cccccc}{\frac13,\frac53,\frac16,\frac56,\frac{5}{12},\frac{5}{12},\frac{-1}{12}}
{1,1,1,~\frac{1}{2},\frac{1}{2},~\frac{3}{2},~\frac{17}{12}\rule[2mm]{0mm}{2mm}}_k
\Big(\frac1{729}\Big)^k,\\
&\Del_k(0)=\frac{7(1+6k)(5+12k)(5616 k^3+11358 k^2+7233 k+1465)}{78732(1+k) (1+2k)^2},\\
&\sig(0)=0
\quad\text{and}\quad
\Gam(\tfrac53,\tfrac53,\tfrac23,\tfrac56)=\frac{15\sqrt3}{8\cdot2^{\frac23}\pi};}
which lead us to the following identity
\[ \frac{59049\sqrt{3}}{14\cdot2^{2/3}\pi}=\sum_{k=0}^{\infty}
\hyp{cccccc}{ \frac13,\frac53,\frac56,\frac76,\frac{5}{12},\frac{-1}{12}}
{ 1,1,~2,~\frac{3}{2},\frac{3}{2},\frac{3}{2} \rule[2mm]{0mm}{2mm}}_k
\frac{\sst 5616 k^3+11358 k^2+7233 k+1465}{729^k}.\]
We claim that the above series is the bisection one of the series below
\[\sum_{k=0}^{\infty}\Lam_{k}
\quad\text{for}\quad
\Lam_{k}:=\hyp{ccc}{\frac23,\frac73,-\frac16}
{1,~1,~2\rule[-2mm]{0mm}{2mm}}_k
\Big(\frac{-1}{27}\Big)^k.\]
This can be justified by computing
\[\Lam_{2k}+\Lam_{2k+1}=\hyp{cccccc}
{\frac13,\frac53,\frac56,\frac76,\frac{5}{12},\frac{-1}{12}}
{1,1,~2,~\frac32,\frac32,\frac32\rule[2mm]{0mm}{2mm}}_k
\frac{\sst 5616 k^3+11358 k^2+7233 k+1465}{1458\cdot729^k}.\]
Therefore, we can evaluate the following simpler series
\[\sum_{k=0}^{\infty}\Lam_{k}
=\sum_{k=0}^{\infty}\Big\{\Lam_{2k}+\Lam_{2k+1}\Big\}
=\frac1{1458}\times\frac{59049\sqrt{3}}{14\cdot2^{2/3}\pi}
=\frac{81\sqrt{3}}{28\cdot2^{2/3}\pi}.\qqed\]

Next, we prove the following elegant formula for a Ramanujan--like series.
\begin{exam}[\fbox{$\veps=0:\frac43,1,1,\frac56$}]\label{pi=4/3}
\[\frac{\pi\Gamma^2 (\frac{1}{3}) }{6\Gamma^2(\frac{5}{6})}
=\sum_{k=0}^{\infty}
\Big(\frac{-1}{27}\Big)^k
\hyp{cccccc}{1,\:\frac{1}2,\:\frac{2}3}
{\frac43,\:\frac43,\:\frac43\rule[2mm]{0mm}{2mm}}_k
\big\{3+7 k\big\}.\]
\end{exam}

\emph{Proof.} \
By specializing the parameters in Theorem~\ref{tre-C}
\[\veps=0\quad\text{and}\quad
\big\{a,b,c,d\big\}=\bigg\{\frac43,1,1,\frac56\bigg\}\]
we can compute
\pnq{&\mc{F}(3k)=\hyp{c}
{\frac12,\frac12,\frac23}
{\\[-3.5mm]\frac32,\frac43,\frac43}_{2k}
\frac{\hh{1,1,\frac12,\frac12,\frac13,\frac56}{k}}
{(3k)!~(\frac32)_{3k}}
=\hyp{cccccc}{1,\frac12,\frac13,\frac14,\frac14,\frac34,\frac56}
{\frac23,\frac23,\frac23,\frac54,\frac76,\frac76,\frac76\rule[2mm]{0mm}{2mm}}_k
\Big(\frac1{729}\Big)^k,\\
&\Del_k(0)=\frac{(1+4k)(4368 k^4+9742 k^3+7799 k^2+2588 k+283)}{72(2+3k)^3},\\
&\sig(0)=0
\quad\text{and}\quad
\Gam(\tfrac43,1,1,\tfrac56)=\frac{\pi\Gam^2(\frac13)}{36\Gam^2(\frac56)};}
which give rise to the following identity
\[\frac{16 \pi  \Gamma (\frac{1}{3})^2}{\Gamma(\frac{5}{6})^2}
=\sum_{k=0}^{\infty}
\hyp{cccccc}{1,\frac12,\frac13,\frac14,\frac34,\frac56}
{\frac53,\frac53,\frac53,\frac76,\frac76,\frac76\rule[2mm]{0mm}{2mm}}_k
\frac{\sst 4368 k^4+9742 k^3+7799 k^2+2588 k+283}{729^k}.\]
For the sequence defined by
\[\Lam_{k}:=(3+7k)\hyp{cccccc}{1,\:\frac{1}2,\:\frac{2}3}
{\frac43,\:\frac43,\:\frac43\rule[2mm]{0mm}{2mm}}_k
\Big(\frac{-1}{27}\Big)^k,\]
it is routine to compute the sum of its two consecutive terms
\[\Lam_{2k}+\Lam_{2k+1}
=\hyp{cccccc}{1,\frac12,\frac13,\frac14,\frac34,\frac56}
{\frac53,\frac53,\frac53,\frac76,\frac76,\frac76\rule[2mm]{0mm}{2mm}}_k
\frac{\sst 4368 k^4+9742 k^3+7799 k^2+2588 k+283}{96\cdot729^k}.\]
Hence the afore--displayed series is equivalent to the following simpler series
\[\sum_{k=0}^{\infty}\Lam_{k}
=\sum_{k=0}^{\infty}\Big\{\Lam_{2k}+\Lam_{2k+1}\Big\}
=\frac1{96}\times\frac{16 \pi  \Gamma (\frac{1}{3})^2}{\Gamma(\frac{5}{6})^2}
=\frac{\pi\Gamma(\frac{1}{3})^2}{6\Gamma(\frac{5}{6})^2}.\]
This completes the proof of the formula given in Example~\ref{pi=4/3}.\qed

By carrying out the same procedure, we shall evaluate further 31 Ramanujan--like
series in closed forms. Compared with the other existing $\pi$-related series
with convergence rate ``$\frac{-1}{27}$" in the literature
(see \cite{kn:chu11mc,kn:chu18e,kn:chu14mc} for example), all the formulae
recorded below are believed to be new, except for Examples~\ref{pi-1A}
and~\ref{pi-1B}. They are divided into four classes according to their
values and exhibited as examples. In each example, the parameter setting
$\boxed{\veps:a,b,c,d}$ and eventual references will be highlighted in the header.
Furthermore, all the formulae are experimentally checked by an appropriately
devised \emph{Mathematica} package in order to ensure the accuracy.

\subsection{Series for $\pi^{2}$}
\begin{exam}[Chu and Zhang~\cito{chu14mc}:
\fbox{$\veps=0:\frac32,1,1,1$}]
\[\frac{\pi^2}{2}
=\sum_{k=0}^{\infty}\Big(\frac{-1}{27}\Big)^k
\hyp{ccc}{1,\:1}
{\frac43,\frac53\rule[2mm]{0mm}{2mm}}_k
\frac{5+7k}{1+2k}.\]
\end{exam}

\begin{exam}[\fbox{$\veps=1:\frac32,1,1,1$}]
\[9\pi^2
=89+\sum_{k=1}^{\infty}\Big(\frac{-1}{27}\Big)^k
\hyp{ccccccccc}{\frac12,\:1,\:3}
{\frac52,\:\frac53,\:\frac73\rule[2mm]{0mm}{2mm}}_k
\big\{17+14k\big\}.\]
\end{exam}

\begin{exam}[\fbox{$\veps=2:\frac52,1,2,1$}]
\[\frac{1575\pi^2}8
=1960-\sum_{k=0}^{\infty}\Big(\frac{-1}{27}\Big)^k
\hyp{ccccccccc}{\frac12,\:1,\:5}
{\frac92,\:\frac73,\:\frac83\rule[2mm]{0mm}{2mm}}_k
\big\{17+7k\big\}.\]
\end{exam}

\begin{exam}[\fbox{$\veps=2:\frac32,1,1,1$}]
\[675\pi^2
=6600+\sum_{k=0}^{\infty}\Big(\frac{-1}{27}\Big)^k
\hyp{ccccccccc}{\frac32,\:1,3}
{\frac72,\frac73,\frac83\rule[2mm]{0mm}{2mm}}_k
\big\{63+59k+14k^2\big\}.\]
\end{exam}

\subsection{Series for $\pi^{2}/\Gam^3$}
\begin{exam}[\fbox{$\veps=0:\frac12,\frac13,\frac13,\frac13$}]\label{pi+2A}
\[\frac{2\pi^2}{\Gam^3(\frac23)}
=\sum_{k=0}^{\infty}\Big(\frac{-1}{27}\Big)^k
\hyp{ccccccccc}{\frac23,\:\frac23,\:\frac16}
{1,\:\frac43,\:\frac76\rule[2mm]{0mm}{2mm}}_k
\big\{8+21k\big\}.\]
\end{exam}

\begin{exam}[\fbox{$\veps=1:\frac12,\frac23,\frac{-1}3,\frac23$}]
\[\frac{2\pi^2}{\Gam^3(\frac13)}
=\sum_{k=0}^{\infty}\Big(\frac{-1}{27}\Big)^k
\hyp{ccccccccc}{\frac{1}3,\:\frac{1}3,-\frac{1}6}
{\!\!1,~\frac23,~\frac{5}6\rule[2mm]{0mm}{2mm}}_k
\big\{1+21 k\big\}.\]
\end{exam}

\begin{exam}[\fbox{$\veps=0:\frac32,\frac23,\frac23,\frac53$}]
\[\frac{45\pi^2}{\Gam^3(\frac13)}
=\sum_{k=0}^{\infty}\Big(\frac{-1}{27}\Big)^k
\hyp{ccccccccc}{\frac{1}3,\:\frac{7}3,-\frac{1}6}
{1,\:\frac53,\:\frac{11}6\rule[2mm]{0mm}{2mm}}_k
\big\{23+42k\big\}.\]
\end{exam}

\begin{exam}[\fbox{$\veps=1:\frac12,\frac13,\frac13,\frac13$}]
\[\frac{5\pi^2}{3\Gam^3(\frac23)}
=\sum_{k=0}^{\infty}\Big(\frac{-1}{27}\Big)^k
\hyp{ccccccccc}{\frac53,-\frac13,-\frac56}
{1,~\frac13,~\frac{7}6\rule[2mm]{0mm}{2mm}}_k
\big\{5-42k\big\}.\]
\end{exam}

\begin{exam}[\fbox{$\veps=0:\frac32,\frac{-1}3,\frac23,\frac53$}]
\[\frac{55\pi^2}{2\Gam^3(\frac13)}
=\sum_{k=0}^{\infty}\Big(\frac{-1}{27}\Big)^k
\hyp{ccccccccc}{\frac{10}3,-\frac{2}3,-\frac{7}6}
{1,~\frac23,~\frac{17}6\rule[2mm]{0mm}{2mm}}_k
\big\{16+21k\big\}.\]
\end{exam}

\begin{exam}[\fbox{$\veps=2:\frac32,\frac13,\frac43,\frac13$}]
\[\frac{91\pi^2}{16\Gam^3(\frac23)}
=\sum_{k=0}^{\infty}\Big(\frac{-1}{27}\Big)^k
\hyp{ccccccccc}{\frac{11}3,-\frac{1}3,-\frac{5}6}
{1,~\frac43,~\frac{19}6\rule[2mm]{0mm}{2mm}}_k
\big\{23+21 k\big\}.\]
\end{exam}

\begin{exam}[\fbox{$\veps=1:\frac32,\frac13,\frac{1}3,\frac43$}]
\[\frac{175\pi^2}{36\Gam^3(\frac23)}
=\sum_{k=0}^{\infty}\Big(\frac{-1}{27}\Big)^k
\hyp{ccccccccc}{\frac53,\:\frac53,-\frac56,-\frac56}
{1,~\frac13,~\frac{7}6,~\frac{13}6\rule[2mm]{0mm}{2mm}}_k
\big\{25+42 k\big\}.\]
\end{exam}

\begin{exam}[\fbox{$\veps=0:\frac52,\frac23,\frac53,\frac53$}]
\[\frac{825\pi^2}{8\Gam^3(\frac13)}
=\sum_{k=0}^{\infty}\Big(\frac{-1}{27}\Big)^k
\hyp{ccccccccc}{\frac{7}3,\:\frac{7}3,-\frac{1}6,-\frac{1}6}
{1,\:\frac53,~\frac{11}6,~\frac{17}6\rule[2mm]{0mm}{2mm}}_k
\big\{53+42 k\big\}.\]
\end{exam}

\begin{exam}[\fbox{$\veps=1:\frac32,\frac23,\frac23,\frac53$}]
\[\frac{60\pi^2}{\Gam^3(\frac13)}
=\sum_{k=0}^{\infty}\Big(\frac{-1}{27}\Big)^k
\hyp{ccccccccc}{\frac{4}3,-\frac{2}3,-\frac{1}6}
{1,~\frac23,~\frac{11}6\rule[2mm]{0mm}{2mm}}_k
\big\{32+111 k+126 k^2\big\}.\]
\end{exam}


\begin{exam}[\fbox{$\veps=0:\frac32,\frac13,\frac43,\frac43$}]
\[\frac{21\pi^2}{\Gam^3(\frac23)}
=\sum_{k=0}^{\infty}\Big(\frac{-1}{27}\Big)^k
\hyp{ccccccccc}{\frac53,-\frac13,\:\frac16}
{1,~\frac43,~\frac{13}6\rule[2mm]{0mm}{2mm}}_k
\big\{83+195k+126k^2\big\}.\]
\end{exam}

\begin{exam}[\fbox{$\veps=1:\frac32,\frac{-1}3,\frac53,\frac23$}]
\[\frac{715\pi^2}{12\Gam^3(\frac13)}
=\sum_{k=0}^{\infty}\Big(\frac{-1}{27}\Big)^k
\hyp{ccccccccc}{\frac{13}3,-\frac{5}3,-\frac{13}6}
{1,~\:\frac23,~\:\frac{17}6\rule[2mm]{0mm}{2mm}}_k
\big\{13+51 k-126 k^2\big\}.\]
\end{exam}

\begin{exam}[\fbox{$\veps=1:\frac32,\frac23,\frac23,\frac23$}]
\[\frac{35\pi^2}{12\Gam^3(\frac13)}
=\sum_{k=0}^{\infty}\Big(\frac{-1}{27}\Big)^k
\hyp{ccccccccc}{\frac{1}3,\frac{7}3,-\frac{7}6}
{1,\frac23,~\frac{11}6\rule[2mm]{0mm}{2mm}}_k
\frac{7-75 k-126 k^2}{(1-6k)(5+6k)}.\]
\end{exam}

\subsection{Series for $\pi^{-1}$}
\begin{exam}[Chu~\cito{chu18e}:
\fbox{$\veps=0:\frac13,\frac13,\frac13,\frac16$}]\label{pi-1A}
\[\frac{9\sqrt3}{2^{\frac43}\pi}
=\sum_{k=0}^{\infty}
\Big(\frac{-1}{27}\Big)^k
\hyp{cccccc}{\frac13,\:\frac23,\:\frac16}
{1,\:1,\:1\rule[2mm]{0mm}{2mm}}_k
\big\{2 + 21k\big\}.\]
\end{exam}

\begin{exam}[Chu~\cito{chu18e}:
\fbox{$\veps=0:\frac23,\frac23,\frac23,\frac56$}]\label{pi-1B}
\[\frac{27\sqrt3}{2^{\frac53}\pi}
=\sum_{k=0}^{\infty}
\Big(\frac{-1}{27}\Big)^k
\hyp{cccccc}{\frac13,\:\frac23,\:\frac56}
{1,\:1,\:1\rule[2mm]{0mm}{2mm}}_k
\big\{5 + 42k\big\}.\]
\end{exam}

\begin{exam}[\fbox{$\veps=1:\frac13,\frac13,\frac13,\frac16$}]
\[\frac{729\sqrt3}{20\sqrt[3]{2}\pi}
=\sum_{k=0}^{\infty}
\Big(\frac{-1}{27}\Big)^k
\hyp{cccccc}{\frac13,\:\frac83,\:\frac16}
{1,\:2,\:2\rule[2mm]{0mm}{2mm}}_k
\big\{16 + 21k\big\}.\]
\end{exam}

\begin{exam}[\fbox{$\veps=2:\frac53,-\frac13,\frac53,\frac56$}]
\[\frac{2673 \sqrt3}{14\sqrt[3]{4}\pi}
=\sum_{k=0}^{\infty}
\Big(\frac{-1}{27}\Big)^k
\hyp{cccccc}{\frac23,\frac{13}3,\frac{5}6,-\frac{7}6}
{1,~1,~3,~\frac{17}6\rule[2mm]{0mm}{2mm}}_k
\big\{65+42 k\big\}.\]
\end{exam}

\begin{exam}[\fbox{$\veps=1:\frac23,\frac23,\frac{-1}6,\frac23$}]\label{pi-1C}
\[\frac{3\sqrt3}{\pi\sqrt[3]4}
=\sum_{k=0}^{\infty}
\Big(\frac{-1}{27}\Big)^k
\hyp{cccccc}{\frac13,-\frac13,-\frac16}
{1,~1,~1\rule[2mm]{0mm}{2mm}}_k
\big\{1-63 k^2\big\}.\]
\end{exam}

\begin{exam}[\fbox{$\veps=0:\frac43,\frac13,\frac43,\frac76$}]
\[\frac{81\sqrt3}{2^{\frac13}\pi}
=\sum_{k=0}^{\infty}
\Big(\frac{-1}{27}\Big)^k
\hyp{cccccc}{\frac13,-\frac13,\:\frac76}
{1,~1,~2\rule[2mm]{0mm}{2mm}}_k
\big\{35+90 k+63 k^2\big\}.\]
\end{exam}

\begin{exam}[\fbox{$\veps=1:\frac23,\frac23,\frac23,\frac56$}]
\[\frac{2187 \sqrt3}{10\sqrt[3]{4}\pi}
=\sum_{k=0}^{\infty}
\Big(\frac{-1}{27}\Big)^k
\hyp{cccccc}{\frac13,\:\frac{2}3,\:\frac{11}6}
{1,~2,~2\rule[2mm]{0mm}{2mm}}_k
\big\{77+144 k+63 k^2\big\}.\]
\end{exam}

\begin{exam}[\fbox{$\veps=1:\frac13,\frac13,\frac13,\frac76$}]
\[\frac{2187\sqrt3}{14\sqrt[3]{2}\pi}
=\sum_{k=0}^{\infty}
\Big(\frac{-1}{27}\Big)^k
\hyp{cccccc}{\frac43,-\frac13,\frac{13}6}
{1,~2,~2\rule[2mm]{0mm}{2mm}}_k
\big\{65+195 k+126 k^2\big\}.\]
\end{exam}

\begin{exam}[\fbox{$\veps=1:\frac43,\frac13,\frac43,\frac76$}]
\[\frac{6561 \sqrt3}{40\sqrt[3]{2}\pi}
=\sum_{k=0}^{\infty}
\Big(\frac{-1}{27}\Big)^k
\hyp{cccccc}{\frac{4}3,\frac{11}3,\frac{1}6}
{1,\:~2,\:~3\rule[2mm]{0mm}{2mm}}_k
\frac{143+285 k+126 k^2}{(1-3k)(2+3k)}.\]
\end{exam}

\begin{exam}[\fbox{$\veps=1:\frac43,\frac13,\frac43,\frac16$}]
\[\frac{2187 \sqrt3}{440\sqrt[3]{2}\pi}
=\sum_{k=0}^{\infty}
\Big(\frac{-1}{27}\Big)^k
\hyp{cccccc}{\frac13,\frac{14}3,-\frac{5}6}
{\!\!1,\:~2,~3\rule[2mm]{0mm}{2mm}}_k
\frac{196+333 k+126 k^2}{(7+6k)(13+6k)}.\]
\end{exam}


\subsection{Series for $\pi$}
\begin{exam}[\fbox{$\veps=1:\frac23,1,1,\frac16$}]
\[\frac{9\pi\Gamma^2 (\frac{2}{3}) }{\Gamma^2(\frac{1}{6})}
=1+\sum_{k=1}^{\infty}
\Big(\frac{-1}{27}\Big)^k
\hyp{cccccc}{3,\:\frac{1}2,-\frac{2}3}
{\frac23,\:\frac{5}3,\:\frac{5}3\rule[2mm]{0mm}{2mm}}_k
\big\{13+21k\big\}.\]
\end{exam}

\begin{exam}[\fbox{$\veps=1:\frac43,1,1,\frac56$}]
\[\frac{4\pi\Gamma^2 (\frac{1}{3}) }{\Gamma^2(\frac{5}{6})}
=71+\sum_{k=1}^{\infty}
\Big(\frac{-1}{27}\Big)^k
\hyp{cccccc}{3,\:\frac{1}2,\:\frac{2}3}
{\frac43,\:\frac{7}3,\:\frac{7}3\rule[2mm]{0mm}{2mm}}_k
\big\{23+21k\big\}.\]
\end{exam}

\begin{exam}[\fbox{$\veps=0:\frac23,1,1,\frac16$}]
\[\frac{3\pi\Gamma^2 (\frac{2}{3}) }{\Gamma^2(\frac{1}{6})}
=\sum_{k=1}^{\infty}
\Big(\frac{-1}{27}\Big)^k
\hyp{cccccc}{1,\:\frac{1}2,-\frac{2}3}
{\frac23,\:\frac{2}3,\:\frac{5}3\rule[2mm]{0mm}{2mm}}_k
k\big\{13+21 k\big\}.\]
\end{exam}

\begin{exam}[\fbox{$\veps=0:\frac13,1,1,\frac{-1}6$}]
\[\frac{\pi\Gamma^2 (\frac{1}{3}) }{12\Gamma^2(\frac{5}{6})}
=\sum_{k=0}^{\infty}
\Big(\frac{-1}{27}\Big)^k
\hyp{cccccc}{1,\:\frac{1}2,-\frac{4}3}
{\frac13,\:\frac13,\:\frac43\rule[2mm]{0mm}{2mm}}_k
\big\{21 k^2+8 k-3\big\}.\]
\end{exam}

\begin{exam}[\fbox{$\veps=0:\frac53,1,1,\frac76$}]\label{pi=5/3}
\[\frac{48\pi\Gamma^2 (\frac{2}{3}) }{\Gamma^2(\frac{1}{6})}
=\sum_{k=0}^{\infty}
\Big(\frac{-1}{27}\Big)^k
\hyp{cccccc}{1,\:\frac{1}2,\:\frac{1}3}
{\frac53,\:\frac{5}3,\:\frac{5}3\rule[2mm]{0mm}{2mm}}_k
\big\{9+28 k+21 k^2\big\}.\]
\end{exam}


\subsection*{Concluding Comments} \
There exist different ways to invert Dougall's $_7F_6$-sum
through \eqref{inv+w} and \eqref{inv-w}. However, all the
dual series that we detected by \emph{Mathematica} are ugly
because of the presence of very complicated weight polynomials.
Here is a couple of discouraging examples.

By examining another triplicate form of Dougall's $_7F_6$-sum
\pnq{\Ome_n\big(a;b+\pavv{\tfrac{1+n}3},c,d+\pavv{\tfrac{1+2n}3}\big)
=&\hyp{c}{1+a-c-d,b+c-a}{1+a-d,b-a}_{\pavv{\frac{1+n}3}}\\
\times\hyp{c}{1+a,\quad b+d-a}{1+a-c,b+c+d-a}_n
&\hyp{c}{1+a-b-c,c+d-a}{1+a-b,d-a}_{\pavv{\frac{1+2n}3}},}
we can arrive, under the parameter setting $\veps=1$ 
and $\{a,b,c,d\}=\Big\{\frac52,2,1,\frac54\Big\}$
and after a long and tedious computations,
at the following series for $\pi$:
\pnq{\frac{75\pi}{8}=30&+\sum_{k=1}^{\infty}\Big(\frac{16}{729}\Big)^k
\hyp{cccccccc}
{1,-\frac12,\frac16,-\frac16,\frac18,-\frac{1}{8},-\frac{3}{8},-\frac{5}{8}\\[-3mm]}
{\frac13,~\frac23,~\frac34,~\frac54,~\frac{7}{12},\frac{11}{12},\frac{13}{12},\frac{17}{12}}_k\\
&\times\Big\{60-101k+1075 k^2-4840 k^3-49360 k^4+136896 k^5\Big\}.}

Analogously, by specifying parameters $\veps=0$ and
$\{a,b,c,d\}=\Big\{\frac12,\frac12,\frac12,\frac13\Big\}$,
we get another series for $\pi^{-1}$:
\pnq{\frac{1485\sqrt3}{\pi}
&=\sum_{k=0}^{\infty}\Big(\frac{16}{729}\Big)^k
\hyp{cccccccc}
{\frac12,\frac23,\frac14,\frac34,~\frac16,~\frac19,\frac49,\frac79\\[-3mm]}
{1,1,1,\frac{4}{3},\frac{4}{3},\frac{17}{18},\frac{23}{18},\frac{29}{18}}_k\\
&\times~\Big\{812 + 20373k + 169774k^2 + 634857k^3 + 1091016k^4 + 693036k^5\Big\}.}


\tableofcontents\end{document}